\documentstyle[12pt]{amsart}

\title{Metric Entropy of Homogeneous Spaces}
\thanks{Partially supported by a grant from the 
National Science Foundation.
Research at MSRI is supported in part by NSF grant DMS-9022140.}

\author{Stanislaw J. Szarek}
\address{\hskip-\parindent
Department of Mathematics \\
Case Western Reserve University \\ Cleveland Ohio 44106}
\email{sjs13@@po.cwru.edu}
\curraddr{Equipe d'Analyse\\
        Universit\'e Paris VI\\
        Tour 46-00, Bte 186\\
        4, Place Jussieu\\
        75252 Paris, FRANCE {\normalfont (through Feb 1997)};
        Mathematical Institute\\
        Polish Academy of Sciences\\
        ul.\ Sniadeckich 8\\
        00-950 Warsaw,   P.O.Box 137\\
        POLAND (Feb--June 1997)}

\begin{document}

\begin{abstract}
For a (compact) subset $K$ of a metric space and $\varepsilon > 0$, the
{\em covering number}
$N(K , \varepsilon )$ is defined as the smallest number of balls of
radius $\varepsilon$ whose union covers $K$.  Knowledge of the {\em metric entropy}, 
i.e., the asymptotic behaviour of
covering numbers for (families of) metric spaces is important in many areas of
mathematics (geometry, functional analysis, probability, coding theory, to name a few).
In this paper we give asymptotically correct estimates for covering numbers
for a large class of
homogeneous spaces of unitary (or orthogonal) groups with respect to some natural metrics,  
most notably the one induced by the operator norm.  
This generalizes earlier author's results concerning covering numbers of
Grassmann manifolds; the generalization is motivated by applications to
noncommutative probability and operator algebras. In the process we give a characterization 
of geodesics in $U(n)$ (or $SO(m)$) for a class of non-Riemannian metric structures.
\end{abstract}

\maketitle

\section{Introduction}

If $(M, \rho )$ is a metric space,  $K \subset M$ a compact subset and $\varepsilon > 0$, the
{\em covering number}
$N(K, \varepsilon ) = N (K, \rho , \varepsilon )$ is defined as the smallest number of balls of
radius $\varepsilon$ whose union covers $K$.   If $K$ is a ball of radius $R$ in a normed space 
of real dimension $d$,  it is easily shown (by a volume comparison argument) that,  
for any $\varepsilon \in (0, R]$,  

\begin{equation}
(R / \varepsilon )^d \leq \ N(K, \varepsilon ) \leq (1 + 2R/\varepsilon )^d.
\end{equation}

The lower and the upper estimate in (1) differ roughly by a factor of $2^d$,  and for 
many applications such an accuracy is sufficient.  On the other hand,  determining more 
precise asymptotics for covering numbers and their "cousins",  {\em packing numbers}  (see section 2), 
e.g. for Euclidean balls is a 
nontrivial proposition (and a major industry).  In this paper we attempt to obtain 
estimates of type (1) for homogenous spaces of the
orthogonal group $SO(n)$ or the unitary group $U(n)$  (like e.g. the Grassmann manifold 
$G_{n,k}$ or the Stieffel manifold, equipped with some natural
metric; we admit metrics induced by unitary ideal norms of matrices,  most notably the 
operator norm).  A typical result will be:  if $M$ is a ``nice" homogenous space of  
$SO(n)$ or  $U(n)$ and $\varepsilon \in (0, \theta(M)]$  (where $\theta(M)$ is some 
{\em computable} ``characteristic" of $M$,  in the more ``regular" cases 
$\theta(M) \approx diam \, M$),  then

\begin{equation}
(c \, diam \, M/ \varepsilon )^d \leq N(M, \varepsilon) \leq (C \, diam\, M/\varepsilon )^d,
\end{equation}

\noindent where $d$ is the (real) dimension of $M$  and  $c$ and $C$ are constants 
independent of 
$\varepsilon $ and  (largely) of $M$.
Of course universality of the constants in question is {\em the} crucial point.

We note in passing that (again,  by a standard ``volume comparison" argument) (2) is 
equivalent to the (normalized)
Haar measure of a ball of radius $\varepsilon$ being between 
$(c_1 \varepsilon /diam \, M)^{\dim M}$ and $(C_1 \varepsilon /diam \, M)^{\dim M}$,  
where $c_1, C_1 > 0$ are {\em universal} constants related to $c, C$.  We also point 
out that the underlying metric being typically non-Riemannian,  the methods of 
Riemannian geometry do not directly apply.

This paper is an elaboration of the note  \cite{Iowa} by the author,  where (2) was proved 
in the special case of $M$
being $SO(n)$,   $U(n)$ or a Grassmann manifold $G_{n,k}.$  (Another argument was 
given later in \cite{pajor}.) The original motivation and application was the 
finite-dimensional basis problem;  more precisely, (2) was used in the proof of 
(Theorem 1.1 in \cite{fdbp}): 

\medskip\noindent {\em There is a constant} $c > 0$
{\em such that,
for every positive integer }$n$, {\em there exists an } $n${\em-dimensional normed space} $B$ 
{\em such that, for every projection} $P$ {\em on }$B$ {\em with, say, } 
$.01 n \leq rank P \leq .99 n${\em , we have }
$\| P:B \rightarrow B \| > c \sqrt{n}$.

\medskip The results of \cite{Iowa} were subsequently applied to other problems in convexity,  
local theory of Banach spaces,  operator theory,  noncommutative probability 
and operator algebras  
(cf. e.g. \cite{herrero}, \cite{exotic}, \cite{cube}, \cite{voiculescu}).   
It turned out recently (see \cite{voiculescu}) that some questions from the last 
two fields lead naturally to queries about validity of 
estimates of type (2)  in settings more general than that of \cite{Iowa}  (which,  
additionally,  had a rather limited circulation).  It is the purpose of this paper to 
provide a {\em reasonably general} answer to such questions.   However,  describing 
asymptotics for {\em completely general}  homogenous spaces of $SO(m)$ or $U(n)$ is,  
in all likelihood,  hopeless.  In this paper we cover a number of special cases, 
including those that 
have been explicitly inquired about.   We identify (easily computable) invariants 
relevant to the problem and provide ``tricks" that could be potentially useful 
to handle cases not addressed here.  Ecxept for brief comments here and there, 
we restrict our attention to $M = G/H$,  
where $G = SO(m)$ or $U(n)$, and $H$ - a connected Lie subgroup of $G$;  but clearly 
most of our  analysis can be extended to other compact linear Lie groups.

The organization of the paper is as follows.  

Section 2 explains notation and presents various preliminary results concerning 
covering numbers and their relatives,  unitarily invariant norms and the exponential map.  
Most of these are known or probably known.  In section 3 we show 
(in Proposition 6) that cosets of 
one parameter semigroups are the geodesics in $U(n)$ (or $SO(m)$) endowed with an 
intrinsic metric induced by a unitarily invariant norm;  this result  
could be of independent interest. In section 4 we discuss several simple examples 
that point out possible obstructions to estimates of type (2) and suggest invariants 
mentioned above. We then use the results of the preceding two sections to show 
estimates of type (2) for an abstract class of homogeneous spaces 
that contains $U(n)$, $SO(n)$ and $G_{n,k}$ (Theorem 8).  In section 5 we discuss various 
possibilities for relaxation of assumptions from section 4,  in particular we  
cover the special cases motivated by applications.  At the end of section 5 
we briefly address the issue of extending the results to metrics generated 
by unitarily invariant norms other than the operator norm.

{\em Acknowledgement.} The author would like to  express his gratitude to D. Voiculescu,  whose 
encouragement was instrumental both in the inception and the completion 
of this work.  The final part of the research has been performed while the author was in 
residence at MSRI Berkeley; thanks are due to the staff of the institute 
and the organizers of the Convex Geometry semester for their hospitality and support.

\section{Notation and Preliminaries}

We start with several remarks clarifying the relationship between covering numbers,
packing numbers and their slightly different versions that exist on the market.  First, since 
the centers of balls in the definition of $N(K,\varepsilon)$, as given in the introduction, do not
necessarily need to be in $K$,  the exact value of $N(K, \cdot)$ may depend on the ambient 
metric space $(M,\rho)$ containing $K$. Accordingly,  it is sometimes more convenient to allow sets of 
{\em diameter} $\leq 2 \varepsilon$ in place of balls 
of {\em radius} $\varepsilon$;  call the resulting the quantity $N'(K,\varepsilon)$.  If the 
centers of the balls in the definition {\em are} required to be in $K$,  call the quantity
$N''(K, \varepsilon)$.  Finally,  let the {\em packing number} $\tilde{N}(K,\varepsilon)$ be 
defined as the maximal cardinality of an 
$\varepsilon$-separated  (i.e.  $\rho(x,x') > \varepsilon$ if $x \neq x'$)
set in $K$.  The quantities $N$, $N'$ and $\tilde{N}$ are related as follows

\begin{equation}
N'(\cdot,\varepsilon) \leq N(\cdot,\varepsilon) \leq N''(\cdot,\varepsilon) 
\leq \tilde{N}(\cdot,\varepsilon) 
\leq N'(\cdot,\varepsilon/2).
\end{equation}

\noindent Consequently,  for our ``asymptotic" results the four quantities 
are essentially interchangeable.

 If the metric space  $(M,\rho)$ is
actually a normed space with a norm $\| \cdot \|$ and unit ball  $B$,
we may write $N(K, \| \cdot \| , \varepsilon )$ or $N(K, B , \varepsilon )$
instead of $N(K, \varepsilon )$ or $N(K, \rho , \varepsilon )$. The technology 
for estimating covering/packing 
numbers of subsets (particularly convex subsets) of normed spaces is quite well-developed 
and frequently rather sophisticated  (see \cite{carl}, \cite{pisier}).  We quote here a 
simple well-known result  
that expresses $N(\cdot , \cdot )$ in terms of a ``volume ratio", and of which
(1) is a special case.

\medskip\noindent {\bf Lemma 1.}  {\em Let } $K$, $B \subset R^d$ {\em with } $B$ 
{\em convex symmetric.  Then, for any } $\varepsilon >0$,
$$( \frac{1}{\varepsilon} )^d \frac{vol \, K}{vol \, B} \leq N(K,B , \varepsilon ) 
\leq N''(K,B , \varepsilon ) \leq N'(K,B,\varepsilon/2)
\leq (\frac{2}{\varepsilon} )^d \frac{vol \,(K + \varepsilon/2 B)}{vol \, B}.$$

The next lemma is just an observation which expresses the fact that the covering/packing 
numbers are invariants of Lipschitz maps.

\medskip\noindent {\bf Lemma 2:} {\em Let } $(M, \rho)$ {\em and } $(M_1, \rho_1)$ 
{\em be metric spaces, } $K \subset M$, $\Phi : K \rightarrow M_1$, {\em  and let }
$L > 0$.  {\em If } $\Phi$ {\em verifies}
$$\rho_1( \Phi (x), \Phi (y)) \leq L \, \rho (x,y) \, \mbox{ {\em for} } \, x,y \in K$$
{\em (i.e. } $\Phi$ {\em is a Lipschitz map with constant } $L$), 
{\em  then, for every } $\varepsilon > 0$,
$$N''(\Phi (K), \rho_1 , L \varepsilon) \leq N''(K, \rho , \varepsilon ).$$
{\em Moreover, } $N''$ {\em can be replaced by } $N'$ {\em or } $\tilde{N}$ 
{\em and,  if } $\Phi$ {\em can be extended to a function on } 
$M$ {\em that is still Lipschitz with constant } $L$, {\em also by } $N$.

\medskip We now turn to our main interest,  the unitary group $U(n)$, the (special) orthogonal 
group $SO(n)$ and their homogeneous spaces.  (As $O(n)$ is geometrically a disjoint
union of two copies of $SO(n)$,  all statements about $SO(n)$ will easily transfer to $O(n)$.)
Throughout the paper we will reserve the letter $G$ to denote,  depending on the context, 
 $U(n)$ or $SO(n)$. Similarly,  we will reserve the letter $\cal{G}$ to denote the 
Lie algebra of $G$,  the space $u(n)$ or $so(n)$ of skew-symmetric matrices.
Since $G$ and $\cal{G}$ are subsets of $M(n)$  (the algebra of $n \times n$ matrices,  
real or complex as appropriate),  they inherit various metric structures from the 
latter.  In this paper we focus on the one induced by the operator norm 
(as an operator on the Euclidean space,  that is),  but will also consider the Schatten 
$C_p$-norms $\|x\|_p = (tr \,|x|^p)^{1/p}$ with the operator norm 
$\| \cdot \|_{op} = \| \cdot \| _{\infty}$ being the limit case.  We will use the same notation
$\|x\|_p$ for the $\ell_p^n$-norm on ${\bf R}^n$ or ${\bf C}^n$,  but this should not 
lead to confusion.  We will also occasionally mention other unitarily invariant norms
(i.e. verifying $\|x\| = \| uxv \|$ if $x \in M(n)$ and $u,v \in G$), each necessarily associated 
with a symmetric norm on ${\bf R}^n$ (which we will also denote by $\| \cdot \|$) via 
$\|x\| = \|(s_k (x))_{k=1}^n \|$,  where $s_1 (x), \dots, s_n (x)$ are singular numbers 
of $x$.  

If $\rho$ is a metric on $G$ and 
$H \subset G$ a closed subgroup,  we consider 
the homogeneous space $M = G/H$  of left cosets of $H$ in $G$ as endowed with 
the canonical quotient metric  $\rho_M (E,F) = \inf \{ \rho (u,v) \, : \, u \in E, v \in F \}$. 
A fundamental example is that of the Grassmann manifold $G_{n,k}$ of $k$-dimensional subspaces 
of ${\bf R}^n$ (resp. ${\bf C}^n)$: the relevant subgroup $H$ of $SO(n)$ 
(resp. $U(n)$) consists of matrices of the form 
\begin{equation}
\left [ 
\begin{array}{cc}
u_1 & 0 \nonumber \\
0 & u_2
\end{array}
\right ]
\end{equation} 
where $u_1 \in SO(k)$ and $u_2 \in SO(n-k)$  (resp. $ U(k), U(n-k)$) and the identfication 
of cosets of $H$ with the subspaces is via $uH \sim uE_k$,  where $E_k$ is the linear span of the 
first $k$ vectors of the standard basis of ${\bf R}^n$ (resp. ${\bf C}^n)$.

It will be usually convenient to consider the {\em intrinsic}
metric,  which we will again call $\rho$, on $M=U(n)$ (or $SO(n)$, or the homogeneous space): 
 $\rho(u, v)$ is the infimum  of
lengths of curves in $M$ connecting $u$ and $v$.  Since this is going to be relevant later, 
we observe that the infimum may be taken over {\em absolutely continuous} curves 
$\gamma \, : \, [a,b] \rightarrow M$ (the infimum is then achieved, while any 
rectifiable curve parametrized by arc length is absolutely continuous).  Then the length
\begin{equation}
\ell(\gamma) = \int_0^1 \| \gamma'(t)\|dt,
\end{equation}
where  $\| \cdot \|$ is the appropriate (unitarily invariant) norm on $M(n)$ if $M=G$  (otherwise 
$\| \cdot \|$ can be interpreted as a quotient norm on the corresponding quotient of the 
relevant Lie algebra;  cf. (14) in section 4 and comments following it). We also point out that 
all the  metrics we consider being bi-invariant,  any curve in $M$ can be lifted 
(by compactness and elementary properties of the $L_1$-norm) 
to a ``transversal" curve in $G$ of the same length. 
The correct abstract framework for these considerations is that of 
Finsler geometry  (see e.g. \cite{busemann},  
but the manifolds we consider being canonically embedded in normed spaces we can afford to be more 
``concrete". 

For future reference we point out that the (operator) norm distance 
and the corresponding intrinsic distance are related as follows
\begin{equation}
\|u-v\|=|1 - e^{i\,\rho(u,v)}|.
\end{equation}
\noindent This follows from Proposition 6 in the next section;    
however, here we just wish to point out that the two metrics 
differ by a factor of $\pi/2$ at the most and this particular fact is implied by 
the more or less obvious inequalities 
$\rho(u,v) \geq \|u-v\| \geq |1 - e^{i\,\rho(u,v)}|$. 
Since,  by definition, the corresponding 
quotient metrics are distances between cosets,  the corresponding two metric 
structure on homogeneous spaces are related in the way analogous to (6). 
Accordingly, estimates of type (2) will 
transfer easily from one metric to the other,  and the choice of the one  
to work with will only be a matter of  convenience and/or elegance. 

Because of the invariance of the metric under the action of $G$  ($=U(n)$ or $SO(n)$), 
one can give estimates for the covering numbers of $M$ (analogous to those of Lemma 1) 
in terms of the Haar measure of balls (cf. the comment following (2)).  However,  since 
the dependence of the measure of a ball on the radius is much less transparent now than 
in the ``linear" case,  such estimates are not necessarily useful.  To overcome this 
difficulty we ``linearize" the problem via the exponential map (composed with the quotient map 
$q:G \rightarrow M$ if necessary) and then use Lemma 2.  Since we operate in the ``classical" 
context,  the exponential map is the standard one
$$\exp x = e^x = \sum_{k =1}^{\infty} \frac{x^k}{k!} \mbox{ for }  x \in M(n),$$
and it will be normally sufficient to consider the restriction of $\exp$ to 
$\cal{G}$ ($= u(n)$ or $so(n)$), the Lie algebra of $G$ ($=U(n)$ or $SO(n)$).  
In order to be able to apply Lemma 2 we must 
``understand"  $\Phi = q \, \circ \, \exp$;  specifically we need to know for which 
$K \subset \cal{G}$ 
we have  $\Phi(K) = M$ (or at least when $\Phi(K)$ is ``large") and for which $K$ 
the restriction $\Phi_{|K}$ (resp. $\Phi_{| \Phi(K)}^{-1}$)  is Lipschitz. 
Concerning the first point,  
it is well known that, in our context, $\exp (\cal{G })$ $= G$.  Moreover, we have

\medskip \noindent {\bf Lemma 3.} {\em Let } 
$K = \{x \in {\cal{G}} \,: \, \|x\|_{\infty} \leq \pi \}$ 
{\em be the ball of radius } $\pi$ {\em in } $\cal{G}$ {\em in the operator norm.  Then }
\newline {\em (a) } $\exp (K) = G$  
\newline {\em (b) } $\exp$ {\em is one-to-one on the interior of } $K$.

\medskip The above is a special case of a more general fact for Lie groups, 
 but in the present setting 
can be seen directly from the fact that every unitary matrix can be diagonalized,  with
the argument for $SO(n)$ being just slightly more complicated.

Lemma 3 asserts that $G$ resembles,  in a sense, a ball in $\cal{G}$.  
However,  for our purposes 
we need more quantitative information about $\exp$, which we collect in the next lemma.

\medskip\noindent {\bf Lemma 4.} {\em For any unitarily invariant norm and the corresponding 
metric on } $G$ {\em  (extrinsic on intrinsic),  the map } $\exp : \cal{G} $ $\rightarrow G$ 
{\em is a contraction. }
\newline {\em On the other hand, let } $\|\cdot\|$ 
{\em be any unitarily invariant norm and set, } for $\theta > 0$,
$$\phi(\theta) = \inf \{\frac{\|e^x - e^y\|}{\|x-y\|} \,: \, x, y \in {\cal{G}}, \, x \neq y, \, 
\|x\|_{\infty} \leq \theta, \, \|y\|_{\infty} \leq \theta \}.$$
{\em Then } $\phi(\theta) > 0$ {\em if } $\theta < \pi$. {\em  Moreover,  if } 
$\theta \in [0,2\pi/3)$, {\em  then}
$$\phi(\theta) \geq \prod_{k=1}^{\infty} {(1 - |1 - e^{i \theta / 2^k}|)}.$$
{\em  In particular } $\phi(\theta) \geq .4$  {\em if } $\theta \leq \pi/4$.

\medskip\noindent {\bf Proof.} The first assertion is classical for the extrinsic (norm) 
metric and hence follows formally 
for the intrinsic metric.  For the other assertions,  we observe first that since 
the derivative of the exponential map at $0$ is the identity,  
\begin{equation}
\lim_{\theta \rightarrow 0^+} \phi(\theta) = 1.
\end{equation}
Let $x, y$ be as in the definition of $\phi(\theta)$.  We have 
\begin{eqnarray*}
e^x - e^y &= &e^{\frac{x}{2}}(e^{\frac{x}{2}}-e^{\frac{y}{2}}) + 
(e^{\frac{x}{2}}-e^{\frac{y}{2}}) e^{\frac{y}{2}} \\
&=&2(e^{\frac{x}{2}}-e^{\frac{y}{2}}) + (e^{\frac{x}{2}}-I)(e^{\frac{x}{2}}-e^{\frac{y}{2}})
+(e^{\frac{x}{2}}-e^{\frac{y}{2}}) (e^{\frac{y}{2}}-I)
\end{eqnarray*}
and so,  by the ideal property of unitarily invariant norms,
\begin{eqnarray*}
\|e^x - e^y\| &= &\|e^{\frac{x}{2}}-e^{\frac{y}{2}}\| (2 - \|e^{\frac{x}{2}}-I\| - \|e^{\frac{y}{2}}-
I\|) \\
&\geq&\phi(\frac{\theta}{2}) \|x-y\| \cdot (1 - |1 - e^{\frac{i \theta}{2}}|).
\end{eqnarray*}
\noindent Iterating and using (7) we obtain the third (and hence the last) 
assertion of the lemma.  For the second assertion 
($\phi(\theta) > 0$ if $\theta < \pi$,  not used in the sequel),   we just briefly sketch 
the argument for $G = U(n)$. 
Let $\theta \in (0,\pi)$ and $\delta > 0$.  Consider first the case of the operator norm. 
 We need to show that if $A, B$ are 
Hermitian with  spectra contained in $[-\theta, \theta]$ and $\|e^{iA} - e^{iB}\| \leq \delta$, 
then $\| A - B \| \leq C(\theta) \delta$,  where $C(\theta)$ depends only on $\theta$ 
(and not on $A, B$ or $n$). By  \cite{bhatia}, Theorem 13.6, the 
eigevalues of $A$ and $B$ (multiplicities counted) are,  ia a certain precise sense, ``close", 
and so,  by perturbation, we may assume that they are identical;  we may also assume that 
all those eigenvalues are inegral multiples of $\delta$. Let $u \in U(n)$ be such that $B = uAu^{-1}$;  
we need to show that $\|e^{iA}u - ue^{iA}\| \leq 4\delta$ implies 
$\|Au - uA\| \leq C(\theta) \delta$ and this follows by writing $u$ as a ``block matrix" in the
spectral subspaces of $A$.  For a general unitarily invariant norm we note that the assertion is 
roughly equivalent to uniform boundedness (with respect to the norm in question and with a bound 
depending only on $\theta$) 
of the inverse of the  derivative of the exponential map.  That derivative is, in a proper 
orthonormal basis,  an antisymmetric ``Schur multiplier" (see \cite{va}, Theorem 2.14.3 
and its proof). 
As a consequence, the inverse is also a ``Schur multiplier" and so its norm with respect to the 
operator norm equals to the norm on the trace class $C_1$ (by duality) and dominates the 
norm with respect to any unitarily invariant norm by interpolation 
(cf. \cite{tomczak}, \S 28 or \cite{gohberg}).

\medskip\noindent {\bf Remark.}  In all likelihood,  a version of Lemma 4 (and of Lemma 5 
that follows) should be known,  at least for the operator norm,  but we couldn't find a reference. 
It would be nice to have an elegant proof which gives good constants in the full range of 
$\theta$ ($\in (0, \pi)$).  We point out that the 
``Schur multiplier" argument indicated above provides a simple ``functional calculus" proof 
(with good constants) in the case of the Hilbert-Schmidt $C_2$-norm.

\medskip\noindent {\bf Lemma 5.} {\em Let } $G = U(n)$ {\em (resp.} $SO(n)$) {\em and } 
$\rho$ {\em the intrinsic metric on } $G$ {\em induced by a unitarily invariant norm } 
$\| \cdot \|$ {\em on } $M(n)$ {\em . Then,  for any } $x, y \in \cal{G}$,
$$\rho(e^{x+y}, e^x e^y) \leq \|[x,y]\|.$$

\medskip\noindent {\bf Proof.}  The argument is similar to,  but slightly more complicated than that of 
the previous lemma.  Denote,  for $t \geq 0$, 
$$\psi(t) = max \{\rho(e^{t(x+y)}, e^{tx} e^{ty}),\rho(e^{t(x+y)}, e^{ty} e^{tx})\}.$$
\noindent Clearly $\psi(0)=0$.  Moreover,  expanding the exponentials and noting that 
$\rho(u,v) / \|u-v\| \rightarrow 1$ as $\|u-v\| \rightarrow 0$  (this follows  
from (6),  but can also be seen from the inequalities in the paragraph following (6), 
which do not depend on Proposition 6) we conclude that 

\begin{equation}
\lim_{t \rightarrow 0} \frac{\psi(t)}{t^2} = \|[x,y]\|;
\end{equation}

\noindent Now
\begin{eqnarray*}
\rho(e^{x+y}, e^x e^y) & \leq \rho(e^{x+y}, e^{\frac{x+y}{2}} e^{\frac{x}{2}} e^{\frac{y}{2}}) + 
\rho(e^{\frac{x+y}{2}} e^{\frac{x}{2}} e^{\frac{y}{2}},
e^{\frac{x}{2}} e^{\frac{y}{2}} e^{\frac{x}{2}} e^{\frac{y}{2}}) \\
&+
\rho(e^{\frac{x}{2}} e^{\frac{y}{2}} e^{\frac{x}{2}} e^{\frac{y}{2}},
e^{\frac{x}{2}} e^{\frac{x+y}{2}} e^{\frac{y}{2}})  +
\rho(e^{\frac{x}{2}} e^{\frac{x+y}{2}} e^{\frac{y}{2}},e^x e^y) \\
& = 3\rho(e^{\frac{x+y}{2}}, e^{\frac{x}{2}} e^{\frac{y}{2}})
+ \rho(e^{\frac{y}{2}} e^{\frac{x}{2}}, e^{\frac{x+y}{2}})
\end{eqnarray*}

\noindent Hence $\psi(1) \leq 4 \psi(\frac{1}{2})$ and,  by the same argument, 
$\psi(t) \leq 4 \psi(\frac{t}{2})$ or $\frac{\psi(t)}{t^2} 
\leq \frac{\psi(\frac{t}{2})}{(\frac{t}{2})^2}$ 
for $t \geq 0$.  In combination with (8) this implies the lemma.

\section{Some Non-Riemannian Geometry}

Our last auxiliary result involves non-Riemannian geometry of $G$ = $U(n)$ (or $SO(n)$).  
It is very well known (in a much more general context) that if $G$ is endowed with 
a bi-invariant Riemannian structure (which is, in our case, the one induced by the 
Hilbert-Schmidt $C_2$-norm on $M(n)$), then the geodesics of $G$ are exactly the cosets of 
one-parameter subgroups (see \cite{helgason}, p. 148, Ex. 5, 6). It is not immediately 
clear how general is this phenomenon. 
Since geodesics are normally defined via affine connections, we should make clear 
here that we emphasize the ``metric" approach: a curve in a manifold $M$ endowed with a metric 
is a geodesic if it {\em locally} realizes the (intrinsic) distance between points as 
explained in the paragraph containing (4)  (the argument can be presumably rewritten by 
starting from the affine connection induced by the group structure,  though). We have

\medskip\noindent {\bf Proposition 6.} {\em Let } $\| \cdot \|$ 
{\em be a unitarily invariant norm on } 
$M(n)$ {\em and } $\rho$ {\em the induced intrinsic metric on } 
$G$ ($ = U(n)$ {\em or } $SO(n)$). {\em Then {\em 
\newline (a) {\em cosets of one parametric semigroups (i.e. curves of the form }
$\gamma (t) = ue^{tx}$, $u \in G$, $x \in \cal{G}$) {\em are geodesics in } $(G, \rho)$
\newline (b) {\em if } $\| \cdot \|$ {\em is strictly convex (which happens in particular if }
$\| \cdot \| = \| \cdot \|_p$ {\em for some } $p \in (1, \infty)${\em ),  
then all geodesics are, up to a change of parameter, of the form given in 
(a) (or arcs of curves of such form) }
\newline (c) {\em if, furthermore, the spectrum of } $u^{-1}v$ 
{\em does not contain } $-1${\em , the curve 
of shortest length (geodesic arc) connecting } $u$ {\em and } $v$ {\em is unique.}

\medskip {\bf Remarks.} (1) A unitarily invariant norm on $M(n)$ is  strictly convex 
(i.e. the unit sphere $\{x \, : \, \| x \| = 1\}$ 
does not contain a segment) iff the associated 
symmetric norm on ${\bf R}^n$ is (cf. \cite{gohberg}). The operator norm and the 
trace class $C_1$-norm are {\em not } strictly convex.
\newline (2) It is a somewhat delicate issue how smooth should be the curves we consider,  
particularly since the results from \cite{bhatia}, to which we refer quite heavily, do not 
fit {\em precisely } our needs {\em exactly} in that respect.  
For the purpose of following the proof below,  the reader should 
think of all functions as being at least $C^1$. As indicated in the previous section,  absolutely 
continuous functions provide a convenient framework;  we comment on the fine points of 
the present argument and on their relevance to \cite{bhatia} at the end of the proof.

\medskip \noindent{\bf Proof.} Since $SO(n)$ is a Lie subgroup of $U(n)$, 
it is enough to prove the lemma for the latter. 
\noindent (a) By unitary invariance of the metric it is enough to prove the assertion for 
``short arcs of one parameter semigroups",  i.e. curves of the form 
$\gamma_0 (t) = e^{tx}$, $t \in [0, 1]$ where $\| x \|$ is ``small".  
By  (5),  the length $\ell (\gamma_0 )$ equals $\| x \|$.  We need to show that any curve $\gamma$
in $U(n)$ connecting $\gamma_0 (0) = I$ and $\gamma_0 (1) = e^{x}$ is of length $\geq \| x \|$. 
The argument is based on  (and in fact very close to) the results on spectral 
variation of unitary matrices presented in \cite{bhatia}, \S 13, 14. Indeed,  Theorem 14.3 and 
Remark 14.4 of \cite{bhatia}, when specified to unitary matrices, say in effect that the map 
$\Sigma$ associating 
to a matrix its spectrum is a contraction when considered as acting from $(U(n), \rho)$ to 
$({\bf C}^n, \| \cdot \|)/S_n$,  where $\| \cdot \|$ is now the symmetric norm 
associated to the unitarily invariant norm in question (a priori defined on ${\bf R}^n$, but 
canonically extendable to ${\bf C}^n$) and $S_n$ is the symmetric group acting on ${\bf C}^n$ by 
permuting the coordinates (we recall that we are working with the {\em intrinsic} metric).  
This in turn implies that, for any curve $\gamma$ in $U(n)$ connecting $I$ and $e^{x}$, we have  
\begin{equation}
\ell (\gamma) \geq \ell (\Sigma(\gamma)).
\end{equation}
  Since $\gamma$ lies in $U(n)$,  
its image  under $\Sigma$ lies in ${\bf T}^n/S_n$,  where 
${\bf T} = \{z \in {\bf C} \,: \, |z| = 1 \}$.  Let now 
$E \,: \, {\bf R}^n \rightarrow {\bf T}^n$ be the exponential map of the group ${\bf T}^n$ 
(i.e., $E((\xi_k)_{k=1}^n) = (\exp(i \xi_k))_{k=1}^n$) considered as acting from 
$({\bf R}^n, \| \cdot \|)$ to ${\bf T}^n$ equipped with the intrinsic metric inherited from 
$({\bf C}^n, \| \cdot \|)$.  Then $E$ is a local isometry,  e.g. it is an isometry when 
restricted to $[-\pi/2, \pi/2]^n$.  Consequently,  the appropriately restricted induced map 
$\tilde{E} \,: \, {\bf R}^n / S_n \rightarrow {\bf T}^n/S_n$ is also an isometry and so, 
if  $\| x \|$ is ``small", 
the length of $\gamma_1 = \tilde{E}^{-1} (\Sigma(\gamma))$ is the same as that of  
$\Sigma(\gamma)$,  in particular,  by (9), $\ell (\gamma_1) \leq \ell (\gamma)$.

Now observe that $\gamma_1$ connects $0$ ($\in {\bf R}^n$)  and 
$\lambda = (\lambda_1, \dots, \lambda_n)$,  where $(i \lambda_1, \dots, i \lambda_n)$ are 
eigenvalues of $x$ with multiplicities.  Clearly  $\| x \| = \| i \lambda \| = \| \lambda \|$ 
and they are all equal to the distance between $0$ and $\lambda$ in ${\bf R}^n/S_n$.  Accordingly
\begin{equation}
\| x \| \leq  \ell (\gamma_1) \mbox{ (in } {\bf R}^n/S_n \mbox{) } 
= \ell (\Sigma(\gamma)) \mbox{ (in } {\bf T}^n/S_n \mbox{) }
\leq  \ell (\gamma) \mbox{ (in } (U(n),\rho) \mbox{)}
\end{equation}
and the assertion (a) is proved.

(b) Again,  it is enough to consider curves of the form 
$\gamma_0 (t) = e^{tx}$, $t \in [0, 1]$ where $\| x \|$ is ``small".  We have to 
show that if $\gamma$ is a curve in $U(n)$ connecting $I$ and $e^{x}$ such
$\ell (\gamma) = \| x \| = \ell (\gamma_0)$,  then (after a change of parameter, if 
necessary), $\gamma = \gamma_0$. In other words,  we need to investigate the cases of 
equality in (10). 
\newline Concerning the first inequality in (10), we observe that since 
$\| \cdot \|$ is strictly convex, a 
straight line segment in $({\bf R}^n,\| \cdot \|)$ with endpoints $\lambda, \lambda'$ 
is {\em strictly} shorter than any other curve connecting these points.  
As a consequence, the same is true 
for the corresponding quotient metric in ${\bf R}^n/S_n$ provided the distance between 
$\lambda, \lambda'$ in the quotient metric equals $\| \lambda - \lambda' \|$ (in particular 
if one of the endpoints is $0$).
Thus we may have equality in the first inequality in (10) only if $\gamma_1$ is a segment,  
i.e. if (after a change of parameter, if 
necessary) the spectrum of $\gamma (t)$, or $\Sigma (\gamma) (t)$, equals 
$(\exp (i \lambda_1 t), \dots, \exp (i \lambda_n t))$,  where,  as before, 
$(i \lambda_1, \dots, i \lambda_n)$ are eigenvalues of $x$.  Let $i  \mu_1, \dots, i \mu_m$ be all 
{\em distinct} eigenvalues of $x$,  and let,  
for $k=1,\dots,m$ and $t \in [0,1]$, 
$P_k (t)$ be the spectral projection of $\gamma (t)$ corresponding to the eigenvalue 
$\nu_k (t) = \exp (i \mu_k t)$  (since $x$ was assumed to be ``small",  
$\nu_k$'s are also distinct). 
If follows from elementary functional calculus that $P_k (t)$ is a 
continuous function of $t$ and,  moreover, has the same smoothness (in $t$) as $\gamma (t)$. 
To prove that $\gamma = \gamma_0$ we need to show that  $P_k (t)$ is constant in $t$ for 
$k=1,\dots,m$.  To this end, we need to analyze the equality case in the second inequality 
in (10), and for that we must go into the proof of Theorem 14.3 of \cite{bhatia}.  In our setting
and notation, it is shown there (p. 68, equation (14.5)) that 
\begin{equation}
\ell (\Sigma(\gamma)) \leq \int_0^1 \|P_{\gamma(t)}\gamma'(t) \| dt 
\leq \int_0^1 \|\gamma'(t) \| dt = \ell (\gamma),
\end{equation}
where,  for $u \in U(n)$,  $P_u$ is a (necessarilly contractive) orthogonal projection in 
$M(n)$ onto the commutant of $u$ (in our case only the restriction of the projection to $\cal{G}$ 
is relevant).  By the strict convexity of $\| \cdot \|$,  we have 
$\|P_u x \| < \| x \|$ unless $x$ is in the 
range of $P_u$  and so the second inequality in (11) 
(and hence the second inequality in (10)) is strict unless $[\gamma'(t),\gamma(t)] = 0$ 
for amost all $t \in [0,1]$.  Now, $\gamma(t)$ being normal, $\gamma'(t)$ must also commute 
with the spectral projections of $\gamma(t)$,  i.e.
\begin{equation}
[\gamma'(t),P_k (t)] = 0  \mbox{ for } k=1, \dots ,m. 
\end{equation}
Since $\gamma(t) = \sum_{j=1}^m \nu_j (t) P_j (t)$,  (11) translates into
\begin{equation}
[\sum_{j=1}^m \nu_j (t) P_j' (t), P_k (t)] = 0  \mbox{ for } k=1, \dots ,m.
\end{equation}
For (small) $h>0$, choose $v(h) \in U(n)$ so that it conjugates 
$\sum_{j=1}^m \nu_j (t) P_j (t)$  and 
$\sum_{j=1}^m \nu_j (t) P_j (t+h)$  (it then also conjugates $P_k (t)$ and $P_k (t+h)$ for
$k=1, \dots ,m$)  
and so that $v(0) = I$ and $v(h)$ depends ``smoothly" on $h$. 
Such a ``smooth" choice is possible wherever $\gamma$ is ``smooth":  as we indicated earlier,  
$(P_k)_{k=1}^m$ is 
then also a ``smooth" curve in the the "generalized Stieffel manifold" (i.e., the quotient 
of $U(n)$ by the commutant of $\gamma(t)$; see section 5 for an 
elaboration of the concept and its framework) and so it can be locally lifted to an 
``equally smooth" curve in $U(n)$.  If $y = y(t) = v'(0)$,  then it is easily seen that  
\begin{equation}
P_k' (t) = [y, P_k (t)]  \mbox{ for } k=1, \dots ,m
\end{equation}
and so (13) translates into 
$$[y, \sum_{j=1}^m \nu_j (t) P_j (t)] \in \mbox{ commutant of } \gamma (t),$$
and so (dropping $t$,  which is fixed, from the formulae), 
$$[y, \sum_{j=1}^m \nu_j P_j] = \sum_{k=1}^m P_k [y, \sum_{j=1}^m \nu_j P_j] P_k ,$$
which is easily seen to be $0$.  On the other hand,  
$$[y, \sum_{j=1}^m \nu_j P_j] = \sum_{k=1}^m \sum_{j=1}^m (\nu_j - \nu_k) P_k y P_j.$$
Since $\{ P_k M(n) P_j,\;  k=1, \dots ,m,\, j=1, \dots ,m \}$ 
form an orthogonal decomposition of $M(n)$, 
the above double sum can be $0$ only if all the terms are $0$.  Recalling that $\nu_k$ are 
distinct,  we deduce that  $j \neq k$ implies $P_k y P_j = 0$.  
Hence $[y, P_k (t)] = \sum_{j \neq k}  P_j y P_k - P_k y P_j = 0$ and so,  by (14), $P_k'(t) = 0$ 
for $k=1, \dots ,m$. 
Since $t$ was an arbitrary point in $[0,1]$ at which $\gamma (t)$ was ``smooth",  it follows that 
$P_k'(t) = 0$ for almost all $t$ and so each $P_k (t)$ is in fact constant,  as required. 
This proves part (b).
\newline (c) This is immediate;  by (b), the geodesic must be (up to a change of parameter) 
of the form $\gamma (t) = ue^{tx},  t \in [0,1]$ with $\ell (\gamma ) = \| x \|$ and 
$e^x = u^{-1}v$.  Since the spectrum of $u^{-1}v$ does not contain $-1$,  the last equation 
has unique solution $x_0$ such that $\| x_0 \|_{\infty} < \pi$,  and any other solution $x_1$ 
verifies $\| x_1 \| > \| x_0 \|$,  thus leading to a {\em strictly } longer curve.

\medskip \noindent{ {\bf Comments.} (1) We need to clarify that the results presented in \cite{bhatia} 
and used above were obtained under ``piecewise $C^1$" assumptions,  both for the curve and 
for (roughly speaking and in our notation) the distance function in 
$({\bf C}^n, \| \cdot \|)/S_n$ applied to $\Sigma(\gamma)$.  Such framework was good enough for 
\cite{bhatia}  (their analysis, even though it fails to identify {\em all} singular points 
of the latter ``distance function" in the general case, can be easily patched by approximating a 
general unitarily invariant norm by a smooth one),  but is insufficient in our context, 
particularly for settling the equality cases. As indicated before,  one solution is to consider 
absolutely continuous curves to insure that variation of a function can be expressed in terms 
of its derivative;  we tacitly used that property several times in our argument).  
Another ``fine point" that makes the argument work is the fact that the ``distance function in 
$({\bf C}^n, \| \cdot \|)/S_n$",  being an minimum of a finite number of convex functions,  
is differentiable almost everywhere (also when composed with an absolutely continuous curve) 
and, moreover,  it  has directional derivatives {\em everywhere}.
\newline (2) In the last part of the argument we did show that,  under certain 
additional assumptions,  if a curve $\gamma$ in $U(n)$ verifies 
$[\gamma(t), \gamma '(t)] = 0$,  then it must be contained in a commutative subgroup.  
It appears  likely that the argument carries over to a more general class of ``sufficiently 
smooth" curves in $U(n)$:  what we seem to be exploiting is that the patern of 
equatities between eigenvalues is constant.

\section{The ``1-Complemented" Subgroups}

Let $G = SO(m)$ or $U(n)$ and $\cal{G}$ ($= so(n)$ or $u(n)$) - the Lie algebra of $G$. 
Let $H$ be a connected Lie subgroup of $G$, $\cal{H}$ - the corresponding 
Lie subalgebra of $\cal{G}$ and  $M = G/H$.  In this 
section we concentrate on the case when $G$ and $M$ are endowed with metric 
structures induced by the operator norm,  we will call the respective 
metrics by $\rho$ and $\rho_M$  (for technical purposes, we may use other unitary 
ideal norms, though).  The purpose of this section is to prove,  in the above context, 
an estimate of type (2) for an abstract class of homogeneous spaces that contains 
$U(n)$, $SO(n)$ and the Grassmannians $G_{n,k}$. 
The argument will depend on a careful analysis of the 
exponential map $\exp \, : \, {\cal{G}} \rightarrow G$ and maps obtained from it;  as an 
illustration we point out here that the following result from \cite{Iowa} is an 
immediate consequence of the results from the preceding two sections.

\medskip\noindent {\bf Theorem 7.} {\em If } $G = SO(n)$ {\em or } $U(n)$ 
{\em (endowed with the operator norm or the induced intrinsic metric } 
$\rho${\em)  and } $\varepsilon \in (0, 2]${\em , then }
$$(\frac {c}{\varepsilon} )^d \leq N(G, \varepsilon) \leq (\frac {C}{\varepsilon} )^d,$$
\noindent {\em where } $d$ {\em is the (real) dimension of }$G$  {\em and } 
$c$ {\em and } $C$ {\em are universal numerical constants.} 

\medskip\noindent {\bf Proof.} By Lemma 3 and the first assertion of Lemma 4,  
$\exp$ is a contractive surjective map 
from the (closed) ball of radius $\pi$ in $\cal{G}$ to $G$.  Consequently, Lemma 2 applied with 
$\Phi = \ exp$ and $L = 1$ and combined with the second inequality in (1) (or,  more precisely, 
the upper stimate on $N''$ given by Lemma 1)
yields the upper estimate for  $N(G, \varepsilon)$.  The lower estimate is obtained 
similarly by applying Lemma 2 to $\Phi = \ exp^{-1}$,  $L = 2.5$ and 
$K = \{u \in G \, : \, \rho(u, I) \leq \pi/4\}$, using the first inequality in (1) 
(or, again more precisely, 
the lower estimate on $\tilde{N}$ given by Lemma 1), 
the last assertion of Lemma 4 and Proposition 6(a).
(Proposition 6(a) is needed only for the definition of $K$ and its use may be avoided here.)

\medskip If $M \neq G$,  the approach will be similar,  but the situation is (necessarily) more 
complicated.  Let $q \,: \, G \rightarrow G/H = M$ be the quotient map and consider 
the short exact sequence 
$0 \rightarrow H \rightarrow G \rightarrow M \rightarrow 0$ and the induced sequence of 
maps between the tangent spaces (at resp. $I \in H$, $I \in G$ and $H \in M=G/H$)
\begin{equation}
0 \rightarrow {\cal{H}} \rightarrow {\cal{G}} \rightarrow T_H M \rightarrow 0
\end{equation}
and so $T_H M$ can be identified with the quotient space $\cal{G}/\cal{H}$;  we mean by that 
{\em isometrically} identified whenever all metric structures are induced by a given 
unitarily invariant norm. Since the derivative of the exponential map at $0$ is the identity 
(in particular an isometry),  we can realize that identification by the canonical 
factorization of the derivative of $q \, \circ \, \exp $ at $0$ (which maps $\cal{G}$ to 
$T_H M$ and vanishes on $\cal{H}$) through $\cal{G}/\cal{H}$.  This shows that (at least small) 
neighbourhoods in $M$ resemble balls in the normed space $\cal{G}/\cal{H}$ and gives 
some heuristic evidence that inequalities of type (2) may hold for $M$.  
However,  for a proof of such an inequality one needs 
``uniform isomorphic" (rather than ``infinitesimal") estimates,  and we will obtain these under 
some additional technical assumptions.  Since the additive structure on $\cal{G}$ and the group 
structure on $G$ are not intertwined by the exponential (or any other) map,  it will be 
more convenient to identify $\cal{G}/\cal{H}$ with $\cal{X} = \cal{H}^{\perp}$ 
(the orthogonal complement of $\cal{H}$ in $\cal{G}$) and to consider 
$\Phi = q \, \circ \, \exp_{|\cal{X}}$,  hoping that the direct sum 
$\cal{X} \oplus \cal{H} = \cal{G}$ is ``well-behaving" with respect to the operator norm 
(or any other unitarily invariant norm that we may need to consider),  which  
happens in many natural examples.  
This leads to our first invariant related to a homogeneous space.  We set

\begin{equation}
\kappa(M) = \| P_{\cal{X}} \| = \| I - P_{\cal{H}} \|,
\end{equation}

\noindent where $P_E$ denotes the orthogonal projection onto $E$ and $\| \cdot \|$ is calculated 
with respect to the operator norm. 

Before stating the results, we will present a simple but illuminating example which 
shows that,  in general,  the ``linearization" of $M$ of the type suggested above may work 
{\em only} on the ``infinitesimal" scale (i.e. only very small neighbourhoods 
are ``equivalent" to balls in the tangent space),
and which leads to one more invariant of $M$.  Let  $G = U(n)$ and $H = SU(n)$.  
It is then easily seen that $M = U(n)/SU(n)$ is {\em isometric } to a circle of radius $1/n$ and 
so covering numbers of $N(M,\varepsilon)$ are ``trivial" if $\varepsilon > \pi/n$. 
(Since $M$ is 1-dimensional,  it is necessarilly ``isotropic" and so there are 
neighbourhoods resembling segments of size comparable to the diameter 
of $M$; in particular (2) still holds. However,  one can also produce ``nonisotropic" 
examples: consider e.g.,  $H = \{ I \} \times SU(n-1) \subset U(n) = G$.)  
The reason for this phenomenon 
is that $SU(n)$ (or ($\cal{H}$,  via the exponential map) is very ``densely woven" 
into $U(n)$.  For example,  $e^{2 \pi i/n}I \in SU(n)$ and 
$\| e^{2 \pi i/n}I - I \| < 2 \pi /n$  (more precisely, 
$\rho (e^{2 \pi i/n}I, I) = 2 \pi /n$ by 
Proposition 6(a)),  even though the shortest path connecting $I$ and $e^{2 \pi i/n}I$ and 
contained in $SU(n)$ is of length $2 \pi (1 -1/n)$ (this follows from the proof of 
Proposition 6(a), the length in question must be 
$\geq$ than the length of the shortest path connecting
$(-2 \pi + 2 \pi/n, 2 \pi/n, \dots,2 \pi/n)$ and $0$ in $\ell_{\infty}^n / S_n$ 
that is contained in the plane $\{ (x_k) \in {\bf R}^n \, : \, \sum x_k = 0 \}$;  
another way to express this is that $e^{2 \pi i/n}I = e^x$ with $x \in \cal{H}$ 
forces $\| x \| \geq 2 \pi (1 -1/n)$). 
To quantify the phenomenon we introduce the following concept.  Given $\theta > 0$,  
we will say that  a closed connected Lie subgroup $H$ of $G = U(n)$ 
(or $SO(n)$) is $\theta$-woven if whenever 
$u \in H$ satisfies $\rho (u,I) \leq \theta$  ($\rho$ is the intrinsic 
metric induced by the operator norm $\| \cdot \|_{\infty}$),  then there exists 
$x \in \cal{H}$, $\| x \|_{\infty} < \pi$ such that $u = e^x$. If $M = G/H$, we set

\begin{equation}
\theta(M) = \sup \{ \theta > 0 \, : \, H \mbox{ is } \theta \mbox{-woven} \} = 
 dist(I,H \backslash \exp (B_{\cal{H}}(\pi)),
\end{equation}

\noindent the distance being calculated using $\rho$.  We then have

\medskip\noindent {\bf Theorem 8.} {\em In the notation above,  assume that } 
$\kappa(M) = 1$. {\em  Then, for any } $\varepsilon \in (0,diam\, M]$,
$$N(M, \varepsilon) \leq (\frac {C diam\, M}{\varepsilon })^d,$$
\noindent {\em where } $d$ {\em is the (real) dimension of } $M$,  $diam\, M$ {\em is 
calculated with respect to } $\rho_M$,  {\em and } $C > 0$ {\em is a universal constant. 
Moreover,  if } $\varepsilon \in (0,\theta(M)/4]$, {\em  then }
$$N(M, \varepsilon) \geq (\frac {c \theta (M)}{\varepsilon })^d ,$$ 
\noindent {\em where } $c > 0$ {\em is a universal constant. The last estimate 
holds also if } $\kappa(M) > 1$, {\em  but the constant $c$ may then depend on } $\kappa(M)$.

\medskip\noindent {\bf Proof.} As suggested earlier,  the proof will involve applying 
Lemma 2 to the (properly restricted) map $q \, \circ \, \exp_{|\cal{X}}$ and its inverse, 
where $\cal{X}$ is the orthogonal complement in $\cal{G}$ of $\cal{H}$  (the Lie 
subalgebra of $\cal{G}$ corresponding to the subgroup $H$),  and will be based on two 
lemmas that follow. Given $r > 0$,  let $B_{\cal{X}}(r)$ be the ball in $\cal{X}$ 
of radius $r$ (with respect to the operator norm) and centered at the origin.  
We then have,  in the notation of Theorem 8:

\medskip\noindent {\bf Lemma 9.}  {\em If } $\kappa(M) = 1$, {\em  then }
$q(exp(B_{\cal{X}}(diam\, M))) = M$. 

\medskip\noindent {\bf Lemma 10.} {\em There exist positive constants } 
$\lambda = \lambda (\kappa(M))$ {\em and } 
$r_0 = r_0(\kappa(M))$ {\em such that if } $r = \min \{r_0,\theta (M)/4 \}$ {\em and }
$x, x' \in B_{\cal{X}}(r)$, {\em  then }
$$\rho_M (q(e^x), q(e^{x'})) \geq \lambda \| x - x' \|.$$

\medskip\noindent {\bf Remarks.} (i) Calculating $\theta (M)$ is not difficult,  
particularly when 
$H$ is semisimple.  Indeed,  suppose $\theta (M) < \pi$  (clearly the maximal 
possible value)
and let $u \in H \backslash \exp (B_{\cal{H}}(\pi))$ be such that $u = e^h$, 
$\| h \|_{\infty} \geq \pi$, while $\| u - I \|_{\infty} = \theta(M) < \pi$, 
in particular $u = e^x$ for some $x \in \cal{G}$, $\| x \|_{\infty} < \pi$.  
Since the commutants of $u$ and $x$ are the same, it follows that $\theta (M)$ is
``witnessed" inside a torus $\cal{T}$ in $G$ containing $u$;  moreover,  $\cal{T}$ 
may be assumed to contain the one-parameter semigroup 
$\{ e^{tx} \, : \, t \in {\bf R} \}$ and to be such that $\cal{T} \cap H$ is 
maximal in $H$.  Consequently,  to determine $\theta (M)$ we only need to examine 
maximal tori in $H$ and their extensions to maximal tori in $G$.  
This is particularly easy if $H$ is semi-simple:  all configurations of the 
tori in question are then 
related by conjugation,  and since the metric we consider is invariant under 
conjugation,  it suffices to check just one such configuration. Such an examination 
will also reveal that $\theta(M) = \pi$,  should that be the case. 
\newline (ii) Since $diam\, U(n) = \pi$ and $q$,  being a quotient map,  
is a contraction,  one {\em always} has $diam\, M \leq \pi$. In any case,  by Proposition 6, 
one can always calculate $diam\, M$ by examining images of one-parameter semigroups 
of $G$ under $q$.
\newline (iii) If $M = G_{n,k}$  (the Grassmann manifold), one verifies directly that 
$\kappa(M) = 1$ and  $diam\, M = \pi/2$.  The former follows from the fact 
that $\cal{X}$ consists of those matrices in $\cal{G}$ ($ = u(n)$ or $so(n)$) 
that are of the form  (cf. (3))
$$\left [ 
\begin{array}{cc}
0 & x  \\
-x^* & 0
\end{array}
\right ] $$
The latter is elementary 
and presumably known:  for two ($k$-dimensional) subspaces $E, F$ of 
${\bf R}^n$ (resp. ${\bf C}^n)$),  $\rho_{G_{n,k}}(E,F)$ is the largest of 
the main angles between $E$ and $F$ (see (4) and comments following it 
for the framework and 
e.g., \cite{Iowa},  p. 174 for the more precise analysis). Finally,  it follows 
immediately from the previous remark that $\theta(G_{n,k}) = \pi$ (the maximal tori in $H$ 
are also maximal in $G$).
\newline (iv) If $\kappa(M) = 1$, one can take in Lemma 10 $r = .12$ and $\lambda = .4$ and if,
moreover, $x'=0$, one can take $r = 5/9$ and $\lambda = .4$;  it follows that 
$q \, \circ \, \exp \,(B_{\cal{X}}(5/9)) \supset \{U \in G/H \, : \, \rho(U, q(I)) \leq 2/9 \}$ 
and that $q \, \circ \, \exp ^ {-1}$ restricted to any of these two sets is Lipschitz 
with constant $4.5$.
\newline (v) The proof gives $\lambda (t)$ and $r_0(t)$ to be of order  $1/t$. 
The argument would be slightly more efficient if we considered $\cal{X}$ as 
endowed with the quotient norm $\cal{G}/ \cal{H}$,  
which is more natural in the context.

Assuming the two lemmas above,  Theorem 8 is shown almost exactly as Theorem 7:  
one applies Lemma 2,  first with $\Phi = q \, \circ \, \exp_{|\cal{X}}$, $L = 1$ and 
$K = B_{\cal{X}}(diam\, M)$ for the upper estimate and then with the inverse 
map restricted to $K = q(exp(B_{\cal{X}}(r)))$ and with $L = \lambda^{-1}$ 
for the lower estimate. (All the fine points are hidden in Lemmas 9 and 10.) 

\medskip\noindent {\bf Proof of Lemma 9.}  We will show that,  for every  $p \in [2, \infty)$,  setting 
$$K_p = \{x \in {\cal{X}} \, :\, \|x\|_p \leq n^{1/p}  diam \, M \}$$
(i.e. $K_p$  is a ball in $\cal{X}$ of radius $n^{1/p} diam \, M$ in the Schatten $C_p$-norm 
$\| \cdot \|_p$), we have 
\begin{equation}
q(\exp (K_p)) = M;
\end{equation}
the assertion of the Lemma will then follow by letting $p \rightarrow \infty$.  To this end,  
observe that since the {\bf R}-linear orthogonal projection $M(n) \rightarrow \cal{G}$ 
is of norm one  (with respect to {\em any} unitarily invariant norm),  $\kappa (M)$ equals 
to the norm of the orthogonal projection from $M(n)$ onto $\cal{X}$.  Now,  since 
$(M(n), \| \cdot \|_p)$ is a complex interpolation space between 
$(M(n), \| \cdot \|_{\infty})$  and $(M(n), \| \cdot \|_2)$,  it follows that $P_{\cal{X}}$ 
is also contractive with respect to the $C_p$-norm  (more generally,  of norm  
$\leq \kappa(M)^{1 - 2/p}$).  Furthermore,  since the  $C_p$-norm is strictly convex 
for $p \in (1, \infty)$,  we conclude that 
\begin{equation}
y \notin {\cal{X}} \Rightarrow \| y - P_{\cal{H}} y \|_p = \| P_{\cal{X}} y \|_p < \| y \|_p .
\end{equation}
For clarity,  we will denote by $M_p$ the manifold $M$ equipped with the quotient metric 
$\rho_{p,M}$ induced by 
the Schatten $C_p$-norm. Note that since the operator norm and the $C_p$-norm 
differ by factor $n^{1/p}$ at the most,  we have $diam\, M_p \leq n^{1/p}diam\, M$.  
Let $gH \in M_p$ and let $\gamma$ be the shortest geodesic in $M_p$ connecting 
$H$ and $gH$,  then $\ell (\gamma) \leq n^{1/p}diam\, M$.  Let $\tilde{\gamma}$ be a 
transversal lifting of $\gamma$ to $G$,  i.e. a
curve in $G$ such that $q \, \circ \, \tilde{\gamma} = \gamma$ and 
$\ell (\tilde{\gamma}) = \ell (\gamma)$.  Then of course $\tilde{\gamma}$ is a geodesic 
in $G$ (with respect to the intrinsic metric $\rho_p$ induced by the $C_p$-norm) 
and without loss of generality we may assume that the initial point of 
$\tilde{\gamma}$ is $I$.  By 
Proposition 6,  $\tilde{\gamma}$ must be (perhaps after a change of parameter) of the form 
$\tilde{\gamma}(t) = e^{ty}, \, 0 \leq t \leq 1$ for some $y \in \cal{G}$ and 
$\ell (\tilde{\gamma}) = \| y \|_p \leq n^{1/p}diam\, M$, and so (18) will follow if we 
show that $y \in \cal{X}$.  Indeed,  if that was not the case, (19) would imply that 
$\| y - P_{\cal{H}} y \|_p < | y \|_p$ and so,  for $t > 0$ sufficiently small we 
would have 
$$\rho_{p,M}(e^{ty} H, H) \leq \rho_p(e^{ty},  e^{tP_{\cal{H}} y}) 
\leq \| ty - tP_{\cal{H}} y \|_p < t \| y \|_p$$
and consequently
$$ \begin{array}{ll}
\ell (\tilde{\gamma}) & = \rho_{p,M}(H,gH) = \rho_{p,M}(H,e^y H) \\
& \leq \rho_{p,M}(H,e^{ty} H) + \rho_{p,M}(e^{ty} H,e^y H) 
< t \| y \|_p + (1-t) \| y \|_p \; ,
\end{array} $$
a contradiction.  This proves Lemma 9.

\medskip\noindent {\bf Proof of Lemma 10.}  We need to show that if 
$\; x, x' \in B_{\cal{X}}(r)$ and $h \in \cal{H}$, then 
$$\Delta \equiv \rho (e^{x'}, e^xe^h) \geq \lambda \| x - x' \|.$$
Since $\| x - x' \| \geq \rho (e^{x'}, e^x) = \rho (e^{-x} e^{x'},I)$ 
and $\Delta = \rho (e^{-x} e^{x'},e^h)$,  it is enough to consider $h \in \cal{H}$ 
such that
$$\rho (e^h,I) \leq (1 + \lambda) \| x - x' \| \leq (1 + \lambda) 2r \leq 4r.$$
If $r \leq \theta (M)/4$  (or just  $2(1 + \lambda)r \leq \theta (M)$),  it follows 
from the definition of $\theta (M)$ (i.e. (17)) that $h \in \cal{H}$ may be 
further assumed to satisfy $\| h \|_{\infty} < \pi$,  hence
$$\| h \|_{\infty} = \rho (e^h,I) \leq (1 + \lambda) \| x - x' \| 
\leq (1 + \lambda) 2r \leq 4r.$$
Now,  by Lemma 4 and Lemma 5,
$$ \begin{array}{ll}
\Delta  &\equiv  \rho (e^{x'}, e^{x}e^h) \geq \rho (e^{x'}e^{-h/2}, e^{x}e^{h/2})  \\
&\geq  \rho (e^{x' - h/2}, e^{x + h/2}) - \| [x', h/2] \| - \| [x, h/2] \| \\
&\geq  \phi (r + \| h \|/2) \| x - x' -h \| - 2r \| h \|\\
&\geq \phi (r + \| h \|/2) \| x - x' \| - 2r \| h \| \\
&\geq  (\phi (r + (1 + \lambda)r) - 2r(1 + \lambda))\| x - x' \| \\
&\geq (\phi (3r) - 4r) \| x - x' \|,
\end{array} $$
\noindent where $\phi ( \cdot )$ is the function from Lemma 4.  
It is now clear from Lemma 4 that if 
$r > 0$ is small enough,  then $\phi (3r) - 4r > 0$.  A more careful 
calculation along the same lines shows that if $r = .12$,  then  $\lambda = .4$ works  
(as indicated in Remark (iv) above).
\newline Finally,  if $\kappa (M) > 1$,  we can only use 
$\| x - x' -h \| \geq \kappa (M)^{-1} \| x - x' \|$ in the fourth inequality in the 
preceding argument.  This results in the last expression being 
$(\phi (3r) \kappa (M)^{-1} - 4r) \| x - x' \|$,  which yields the assertion 
of the Lemma with $r_0$ and $\lambda$ being of order $\kappa (M)^{-1}$.

\section{Extensions and Other Tricks}

The scheme presented in the preceding section yields resonable estimates for 
covering numbers $N(M, \varepsilon)$ (with respect to the metric induced by 
the operator norm) of a homogeneous space $M = G/M$  whenever
$\varepsilon \leq diam\, M$ or $\varepsilon \leq \theta(M)$  (for the upper 
and lower estimate respectively) and whenever $\kappa(M)$ is 
appropriately controlled.  This leaves several cases and gray areas that 
are not covered.
\newline (i) The range $\theta(M) < \varepsilon < diam\, M$ even if $\kappa(M) = 1$ and,  
in general, a clarification of the role of the ratio $diam\, M / \theta(M)$  
(the lower and upper estimates differing roughly by $(diam\, M / \theta(M))^d$).
\newline (ii) The upper estimate whenever $\kappa(M) > 1$,  but still "under control".
\newline (iii) The case when we do not control $\kappa(M)$.

With regards to (i), a modification of the example that led to the definition of $\theta(M)$ 
suggested there ($H = \{ I \} \times SU(n-1) \subset U(n) = G$) shows that it is 
possible for $diam\, M$ and $\theta(M)$ to differ by a large factor (of order $n$ 
in that case).  Even though an analysis of such cases is imaginably possible,  
it would be clearly combinatorial and/or algebraic in nature and we do not 
attempt it here.

Concerning (iii),  it is also conceivable that the phenomenon of having $\kappa(M)$ 
``large" can be ``dissected" and expressed in terms of combinatorial/algebraic 
invariants suggested above,  but,  again,  in the examples motivating this work 
(see below) we have $\| P_{\cal{H}} \| = 1$ and hence $\kappa(M) \leq 2$.

This leaves the gap related to (ii):  the examples with, say,  $1 < \kappa(M) \leq 2$ 
do naturally occur and it would be nice to have, at least for that case, an upper 
estimate for covering numbers of $M$ of the type $(C diam \, M/ \varepsilon)^d$. 
Unfortunately,  we do not know how to settle that question in full generality. 
Instead,  we present a ``bag of tricks" that allow to handle various special cases. 
This,and some comments concering covering numbers relative to metrics generated 
by unitarily invariant norms other than the operator norm constitutes this section.

The first observation is that trying to mimmick the proof of Lemma 9 in the case when 
$\kappa(M) > 1$ one arrives at the following picture. Let 
$Q : \cal{G} \rightarrow \cal{G}/ \cal{H}$ be the quotient map, and consider the 
semi-norm $p$ on $\cal{G}$ defined by $p(x) = \| Qx \|_{\infty}$. 
Let $\Lambda : \cal{G}/ \cal{H} \rightarrow \cal{G}$ be 
a norm-preserving lifting of $Q$ (in general nolinear). The argument mimmicking the 
proof of Lemma 9 connects then geodesics in $M$ with ``rays" in the range of $\Lambda$ 
and we could give upper estimates for entropy of $M$ if we were able to control 
entropy of the range of $M$ (e.g. with respect to the semi-norm $p$).

The two specific subgroups of $U(n)$, for which estimates for covering numbers of the 
respective homogeneous spaces are of interest from the point of view of free 
probability  (cf. \cite{voiculescu}, Remark 7.2 and \cite{voiculescupc}), consist of unitaries of some 
$C^*$-subalgebras of $M(n)$,  namely

\noindent (1) The ``block-diagonal" algebra:  the commutant of $\{ P_1, P_2, \dots, P_m \}$,
where $P_j$'s are orthogonal projections
whose ranges form an orthogonal decomposition of ${\bf C}^n$.
\newline (2) The ``tensor-factor" algebra:  if  $n = mk$,  identify ${\bf C}^n$ with 
${\bf C}^m \otimes{\bf C}^k$ and consider matrices of the form $I \otimes x$, $x \in M(k)$; 
these can be also thought of as block matrices with $m$ identical $k \times k$ blocks 
along the diagonal. 

Let us comment here that the homogeneous space obtained from (1) is a 
``generalized Stieffel manifold" of ``orthogonal frames" of subspaces of 
${\bf C}^n$ (or ${\bf R}^n$) with given pattern of dimensions.
 In both cases (1) and (2) the subgroup $H$ (resp. the Lie algebra $\cal{H}$ consists of 
(all) unitaries 
(resp. skew-symmetric matrices having form (1) or (2)),  and the conditional 
expectation is a norm one (with respect  to any unitarily invariant norm) 
projection from $\cal{G}$,  the Lie algebra of $G = U(n)$, onto $\cal{H}$,  in 
particular $\kappa(G/H) \leq 2$.  However,  except for the case $m = 2$ in 
(1) (i.e. the Grassmann manifold) or (2),  we have $\kappa(G/H) > 1$.

Concerning the other parameters,  it is easily seen that in all cases 
$\theta(G/H) = 2$ and $\pi/2 \leq diam \, G/H \leq \pi$.  Accordingly,  Theorem 8 
gives good lower estimates for the covering numbers of $G/H$ and it remains to 
handle the upper ones.  We will use the following (ad hoc)

\medskip\noindent {\bf Theorem 11.} {\em Let } $\alpha \in (0, 1/2]$;  
{\em let } $n$, $G$ ($= U(n)$ or $SO(n)$), 
$H$, $M = G/H$, $\cal{G}$, $\cal{H}$, {\em and } $d$ {\em be as before and assume that } 
$$\min \{ \theta (M), diam\, M, \kappa(M)^{-1} \} \geq \alpha.$$  
{\em Furthermore,  asuume that one of the following holds } 
\newline (a) $\dim\, H \leq (1 - \alpha) \dim\, G$
\newline (b) $H$ {\em acts reducibly on } ${\bf C}^n$ {\em (resp. } ${\bf R}^n$  
{\em and there is a reducing subspace } $E$ {\em with } 
$\alpha n \leq \dim\, E \leq (1 - \alpha) n$
\newline (c) $H$ {\em acts reducibly on } ${\bf C}^n$ {\em (resp. } ${\bf R}^n)$  
{\em and there is a reducing subspace } $E$ {\em with } 
$\dim\, E \equiv k \geq \alpha n$ {\em and such that the orthogonal decomposition } 
${\bf C}^n = E \oplus E^{\perp}$ {\em induces an isomorphism } 
$H \rightarrow U(k) \times H_0$ {\em for some subgroup } $H_0$ {\em of } 
$U(n-k)$ {\em (resp. } $ {\bf R}^n , SO(k), SO(n-k))$. 
\newline {\em Then, for any } $\varepsilon \in (0,diam\, M]$,
$$(\frac {c}{\varepsilon} )^d \leq N(M, \varepsilon) \leq (\frac {C}{\varepsilon} )^d,$$
\noindent {\em where } $c, C > 0$ {\em are constants depending only on } $\alpha$.

\medskip\noindent {\bf Corollary 12.} {\em If } $H \subset U(n)$ {\em is the group of 
unitaries of a ``block-diagonal" or ``tensor-factor" algebra (described in (1) or (2)),  
then the assertion of Theorem 11 holds (e.g., with constants corresponding to } 
$\alpha = 1/3$ {\em ).} 

\medskip\noindent {\bf Proof.} We may of course assume $m \geq 2$. In the case (2) the 
condition (b) of Theorem 11 is always satisfied (with $\alpha = 1/3$). The same is
true in the case (1) except if one of the projections $P_j$ is of rank $\geq n/3$, 
in which case (c) holds.

\medskip\noindent {\bf Remarks.} (i) As was pointed out in \cite{voiculescu}, 
Remark 7.2, the estimates for covering numbers  given by our Corollary 12 
(the ``block-diagonal" case) 
allow sharp free entropy and free entropy dimension estimates. Similarly,  the
 ``tensor-factor" case of the Corollary implies estimates for free entropy 
and free entropy dimension of certain generators of free product von Neumann 
algebras (\cite{voiculescupc}).
\newline (ii) The ``block-diagonal" case of Corollary 12 implies estimates for 
covering numbers of some sets of matrices needed in \cite{dykema}.

\medskip\noindent {\bf Proof of Theorem 11.} As observed earlier,  it is enough to 
show the upper estimate.
\newline (a) The condition (a) is equivalent to $\dim\, M \geq \alpha \dim\, G$. 
It follows from Theorem 7 that, for any $\varepsilon \in (0, 2]$,  $G$ (hence, 
by Lemma 2, $M = G/H$) admits an $\varepsilon$-net of cardinality 
$ \leq (\frac {C}{\varepsilon} )^{dim \, G}$ .  If  $\varepsilon \geq \beta$ and (a) 
holds, this does not exceed $C(\alpha, \beta)^{\dim\, M}$  (the {\em real} 
dimension).  It follows now from Lemma 10 
that the image of $B_{\cal{X}}(r)$ (where $r = r(\alpha)$) contains a ball in $M$ 
of radius $r_1 = r_1(\alpha)$, and so the former and the latter admit, for 
$\varepsilon \leq r$,  an $\varepsilon$-net of cardinality 
$ \leq (\frac {C_1 r}{\varepsilon} )^{dim \, M}$.  Combining this with the 
preceding observation (applied to $\beta = r_1$) we get the required upper estimate. 
\newline (b) The condition (b)  implies (a)  (with $\alpha (1-\alpha)$ in place 
of $\alpha$).
\newline (c) Since the arguments in the real and complex case are identical, 
we restrict the discussion to the latter. 
The condition (c) is not included in (b) only if 
$\dim\, E \equiv k > (1 - \alpha) n$, in particular $k > n/2$. 
Let $H_1 \subset G$ be the subgroup of the form $U(k) \times \{ I \}$ 
in the sense indicated in the condition (c),  then $H_1 \subset H$ and so 
$M = G/H$ is a quotient of $M_1 = G/H_1$.  Now $M_1$ is isomorphic to 
$G_{n,k} \times U(n-k)$ (for sure Lipschitz isomorphic with constant 2 if 
the product metric on the latter is defined in the ``$\ell_{\infty}$ sense"),  
and so, by Theorems 7 and 8, it admits,  for any $\varepsilon \in (0, 2]$, 
an $\varepsilon$-net of cardinality 
$\leq (\frac {C_2}{\varepsilon} )^{dim \, M_1}$ 
(again, the {\em real} dimension). Since $k > n/2$ implies 
$dim \, M_1 < 2 dim \, M$,  arguing as in (a) we obtain the assertion.

In some applications (see e.g. \cite{exotic}, \cite{cube}) it is important 
to know the metric entropy of $M$ equipped with a metric induced by 
unitarily invariant norms other than the operator norm. The scheme presented 
in this paper can be adapted to yield fairly sharp results in the general case. 
Indeed, Lemmas 4, 5 and 6 involve statements about generic unitarily 
invariant norms. Similarly, Lemma 10 and its proof carry over almost 
word by word to the case of an arbitrary unitarily invariant norm 
$\| \cdot \|$ once the parameters such as $\theta$ and $\kappa$ are 
properly interpreted: the balls $B_{\cal{H}}(\cdot)$, $B_{\cal{X}}(\cdot)$ 
are to remain to be defined by the operator norm $\| \cdot \|_{\infty}$, 
but $\theta (M)$ has to be the distance between $I$ and 
$H \backslash \exp (B_{\cal{H}}(\pi))$ in the intrinsic 
metric on $G$ induced by $\| \cdot \|$;  $\kappa(M)$ may be 
calculated using $\| \cdot \|$ (which results in a quantity not larger 
than the one given by the operator norm,  in particular $\kappa(M)=1$ 
if $\| \cdot \|$ is the Hilbert Schmidt norm). The ``linearization" 
procedure can be then implemented and the problem is reduced to 
estimating covering numbers of balls in $\cal{X}$ in the operator norm 
with respect to $\| \cdot \|$.  As indicated in section 2 (see Lemma 1 and 
the paragraph preceding it),  there exist numerous tools for obtaining 
such estimates. In particular,  in many natural cases (e.g. $M = U(n)$, $SO(n)$ or 
the Grassmann manifolds $G_{n,k}$),  the volumes of a ball in $X$ with respect 
to a unitarily invariant norm $\| \cdot \|$ and the inscribed operator norm ball differ 
by a factor $C^d$, $C$ - an universal constant  (this is easily implied e.g., by classical 
facts from \cite{Santalo}, \cite{chevet}, cf. \cite{fdbp}, p. 162; 
see also \cite{gluskin} or \cite{SR}),  which allows to use Lemma 1 
to show that $\log N(M,\varepsilon) \approx d \, \log(diam \, M /\varepsilon)$  (and 
$diam \, M \approx \| I \|$ if $M = U(n)$ or $SO(n)$,  
$diam \, M \approx \| P \|$,  where $P$ is an orthogonal projection of rank 
equal to $\min \{ k, n-k \}$ if $M = G_{n,k}$);  these cases have been worked out in 
\cite{Iowa}.

\end{document}